\documentclass[12TP,draft]{article}

\begin{document}

\begin{center}
\LARGE\noindent\textbf{On pre-Hamiltonian cycles in balanced bipartite digraphs }\\

\end{center}
\begin{center}
\noindent\textbf{Samvel Kh. Darbinyan }\\

Institute for Informatics and Automation Problems, Armenian National Academy of Sciences

E-mails: samdarbin@ipia.sci.am\\
\end{center}

\textbf{Abstract}

Let $D$ be a strongly connected balanced bipartite directed graph of order $2a\geq 10$. Let $x,y$ be distinct vertices in $D$. $\{x,y\}$ dominates a vertex $z$ if $x\rightarrow z$ and $y\rightarrow z$; in this case, we call the pair $\{x,y\}$ dominating.  In this paper we prove: 

 {\it If the underlying undirected graph of $D$ is not 2-connected and $max\{d(x), d(y)\}\geq 2a-2$ for every dominating pair of vertices $\{x,y\}$,  then  
 $D$ contains a cycle of length $2a-2$ unless  $D$ is isomorphic to a certain digraph of order ten which we specify.}\\


\section {Introduction} 

We consider directed graphs (digraphs) in the sense of \cite{[1]}. A cycle of a digraph $D$ is called Hamiltonian if it contains all the vertices of $D$. For convenience of the reader terminology and notations will be given in details in section 2. A digraph $D$ of order $n$ is Hamiltonian  if it contains a  Hamiltonian cycle and pancyclic if it contains cycles of every length $k$, $3\leq k\leq n$.  For general digraphs there are several sufficient conditions for existence of Hamiltonian cycles in digraphs. In this paper, we will be concerned with the degree conditions.

The well-known and classical are  Ghouila-Houri's, Nash-Williams', Woodall's, Meyniel's and Thomas\\sen's theorems 
(see, e.g., \cite{[2]}- \cite{[6]}). There are analogies results of the above-mentioned theorems for the pancyclicity of digraphs (see, e.g., [7-12]).

 Each of theorems (\cite{[2]}-\cite{[6]}) imposes a degree condition on all pairs of nonadjacent vertices (or on all vertices).
In \cite{[13]} and \cite{[14]},  some sufficient conditions were described for a digraph to be Hamiltonian, in which a  degree condition is required only for some pairs of nonadjacent vertices. Let us recall  only the following theorem of them. \\

 \textbf{Theorem 1.1} (Bang-Jensen, Gutin, Li \cite{[13]}). {\it Let $D$ be a strongly connected digraph of order $n\geq 2$. Suppose that $min\{d(x),d(y)\}\geq n-1$ and  $d(x)+d(y)\geq 2n-1$ for any pair of nonadjacent vertices $x,y$ with a common in-neighbor. Then $D$ is Hamiltonian.}\\

A digraph $D$ is called  bipartite  if there exists a partition $X$, $Y$ of its vertex set into two partite sets such that every arc of $D$ has its end-vertices in different partite sets. 
It is called balanced if $|X|=|Y|$.

A cycle of a non-bipartite digraph $D$ is called pre-Hamiltohian if it contains all the vertices of $D$ except one. 
The concept of pre-Hamiltonian cycle for the balanced bipartite digraphs is the following:

A cycle of a balanced bipartite digraph $D$ is called pre-Hamiltonian if it contains all the vertices of $D$ except two.

There are results analogies to the theorems of Ghouila-Houri, Nash-Williams, Woodall, Meyniel and Thomassen 
for balanced bipartite digraphs (see e.g., \cite{[15]} and the papers cited there).

Let $x,y$ be a pair of distinct vertices in a digraph $D$. We call the pair $\{x,y\}$ dominating, if there is a vertex $z$ in $D$ such that $x\rightarrow z$ and $y\rightarrow z$.

An analogue of Theorem 1.1  for bipartite digraphs was given by R. Wang \cite{[16]} and recently strengthened by the author   \cite{[17]}.

\textbf{Theorem 1.2} (R. Wang \cite{[16]}). {\it Let $D$ be a strongly connected balanced bipartite digraph of order $2a$, where $a\geq 1$. Suppose that, for every dominating pair of vertices $\{x,y\}$, either $d(x)\geq 2a-1$ and $d(y)\geq a+1$ or $d(y)\geq 2a-1$ and $d(x)\geq a+1$. Then $D$ is Hamiltonian}.\\

Let $D$ be a balanced bipartite digraph of order $2a\geq 4$. For integer $k\geq 0$, we say that
 $D$ satisfies condition $B_k$ when $max\{d(x),d(y)\}\geq 2a-2+k$ for every pair of dominating vertices  $x$ and $y$.\\

\textbf{Theorem 1.3} (Darbinyan \cite{[17]}). {\it Let $D$ be a strongly connected balanced bipartite digraph of order $2a$, where $a\geq 4$. Suppose that $D$ satisfies condition $B_1$, i.e., for every dominating pair of vertices $\{x,y\}$, either $d(x)\geq 2a-1$  or $d(y)\geq 2a-1$. Then either $D$ is Hamiltonian or isomorphic to the digraph $D(8)$ (for the definition of $D(8)$, see Example 2)}.\\

A balanced bipartite digraph of order $2m$ is even pancyclic  if it contains a cycle of length $2k$ 
for any $2\leq k\leq m$. 

An even pancyclic version of Theorem 1.3 was proved in \cite{[18]}.\\

\textbf{Theorem 1.4.} (Darbinyan \cite{[18]}). {\it Let $D$ be a strongly connected balanced bipartite digraph of order $2a\geq 8$ other than the directed cycle of length $2a$. If $D$ satisfies condition $B_1$, i.e., $ max\{d(x), d(y)\}\geq 2a-1$ for every dominating pair of vertices $\{x,y\}$,  then either $D$ contains  cycles of all even lengths less than or equal to $2a$ or $D$ is isomorphic to digraph $D(8)$.}\\

\textbf{Theorem 1.5.} (Darbinyan \cite{[18]}). {\it Let $D$ be a strongly connected balanced bipartite digraph of order $2a\geq 8$ which contains a pre-Hamiltonian cycle (i.e., a cycle of length $2a-2$). If $D$ satisfies condition $B_0$, i.e., $max \{d(x), d(y)\}\geq 2a-2$ for every dominating  pair of vertices $\{x,y\}$, then for any $k$, $1\leq k\leq a-1$, $D$ contains  cycles of every length $2k$, $1\leq k\leq a-1$.}\\

In view of Theorem 1.5 it seems quite natural to ask whether a balanced bipartite digraph of order $2a$  in which 
$max \{d(x), d(y)\}\geq 2a-2$ for every dominating pair of vertices $\{x,y\}$  contains a pre-Hamiltonian 
cycle (i.e., a cycle of length $2a-2$). 

The underlying undirected graph of a digraph $D$ is the unique graph such that it contains an edge $xy$ if $x\rightarrow y$ or $y\rightarrow x$ (or both).

In this paper we prove the following theorem.\\

\textbf{Theorem 1.6}. {\it Let $D$ be a strongly connected balanced bipartite digraph of order $2a\geq 10$  with partite sets $X$ and $Y$. Assume that the underlying undirected graph of $D$  is not 2-connected and $D$ satisfies condition $B_0$. Then  $D$ contains a cycle of length $2a-2$ unless $D$ is isomorphic to the digraph $D(10)$ (for the definition of $D(10)$, see Example 2)}.\\

\section {Terminology and Notations}

Terminology  and notations not described below follow \cite{[1]}. 
  In this paper we consider finite digraphs without loops and multiple arcs. The vertex set and the arc set of a digraph $D$  denoted
  by $V(D)$  and   $A(D)$, respectively. The order of $D$ is the number
  of its vertices.  
For any $x,y\in V(D)$, we also write $x\rightarrow y$ if $xy\in A(D)$. If 
$xy\in A(D)$, then we say that $y$ is an out-neighbour of $x$ and $x$ is an in-neighbour of $y$.   The notation $x\leftrightarrow y$ means that $x\rightarrow y$ and $y\rightarrow x$ ($x\leftrightarrow y$ is called 2-cycle). We denote by $a(x,y)$ the number of arcs with end-vertices $x$ and $y$. For disjoint subsets $A$ and  $B$ of $V(D)$  we define $A(A\rightarrow B)$ \,
   as the set $\{xy\in A(D) : x\in A, y\in B\}$ and $A(A,B)=A(A\rightarrow B)\cup A(B\rightarrow A)$. If $x\in V(D)$
   and $A=\{x\}$ we sometimes write $x$ instead of $\{x\}$. If $A$ and $B$ are two disjoint subsets of $V(D)$ such that every
   vertex of $A$ dominates every vertex of $B$, then we say that $A$ dominates $B$, denoted by $A\rightarrow B$. The notation $A\leftrightarrow B$ means that $A\rightarrow B$ and $B\rightarrow A$.

 We let $N^+(x)$, $N^-(x)$ denote the set of  out-neighbours, respectively the set  of in-neighbours of a vertex $x$ in a digraph $D$.  If $A\subseteq V(D)$, then $N^+(x,A)=A\cap N^+(x)$, $N^-(x,A)=A\cap N^-(x)$ and $N^+(A)=\cup_{x\in A}N^+(x)$. The out-degree of $x$ is $d^+(x)=|N^+(x)|$ and $d^-(x)=|N^-(x)|$ is the in-degree of $x$. Similarly, $d^+(x,A)=|N^+(x,A)|$ and $d^-(x,A)=|N^-(x,A)|$. The degree of the vertex $x$ in $D$ is defined as $d(x)=d^+(x)+d^-(x)$ (similarly, $d(x,A)=d^+(x,A)+d^-(x,A)$). The subdigraph of $D$ induced by a subset $A$ of $V(D)$ is denoted by $D\langle A\rangle$ or $\langle A\rangle$ for brevity.   

 The path (respectively, the cycle) consisting of the distinct vertices $x_1,x_2,\ldots ,x_m$ ( $m\geq 2 $) and the arcs $x_ix_{i+1}$, $i\in [1,m-1]$  (respectively, $x_ix_{i+1}$, $i\in [1,m-1]$, and $x_mx_1$), is denoted by  $x_1x_2\cdots x_m$ (respectively, $x_1x_2\cdots x_mx_1$). 
We say that $x_1x_2\cdots x_m$ is a path from $x_1$ to $x_m$ or is an $(x_1,x_m)$-path.  A cycle that contains  all the vertices of $D$  is a Hamiltonian cycle. 
 A digraph $D$ is strongly connected (or, just, strong) if  there exists an $(x,y)$-path in $D$   for every ordered pair of distinct vertices $x,y$ of $D$.
    
Two distinct vertices $x$ and $y$ are adjacent if $xy\in A(D)$ or $yx\in A(D) $ (or both). For integers $a$ and $b$, $a\leq b$, let $[a,b]$  denote the set of
all integers which are not less than $a$ and are not greater than
$b$.

\section {Examples}

\textbf{Example 1.} Let $D(10)$ be a bipartite digraph with partite sets $X=\{x_0,x_1,x_2,x_3,x_4\}$ and 
$Y=\{y_0,y_1,y_2,y_3,y_4\}$ satisfying the following conditions: The induced subdigraph $\langle\{x_1,x_2,x_3,y_0,y_1\}\rangle$ is a complete bipartite digraph with partite sets  $\{x_1,x_2,x_3\}$ and $\{y_0,y_1\}$;  $\{x_1,x_2,x_3\}\rightarrow\{y_2,y_3,y_4\}$; $x_4\leftrightarrow y_1$; $x_0\leftrightarrow y_0$ and $x_i\leftrightarrow y_{i+1}$ for all $i\in [1,3]$. $D(10)$ contains no other arcs.

It is easy to check that the digraph $D(10)$ is strongly connected  and satisfies condition $B_0$, but the underlying undirected graph of $D(10)$ is not 2-connected and $D(10)$ has no cycle of length 8. (It follows from the facts that $d(x_0)=d(x_4)=2$ and $x_0$ ($x_4$) is on 2-cycle). It is not difficult to check that any digraph obtained from $D(10)$ by adding a new arc whose one end-vertex is $x_0$ or $x_4$ contains no cycle of length eight. Moreover, if to $A(D)$ we add some new arcs of the type $y_ix_j$, where $i\in [2,4]$ and $j\in [1,3]$, then always we obtain a digraph which does not satisfy condition $B_0$. \\ 

\textbf{Example 2.} Let  $K_{2,3}^*$ be a complete bipartite digraph with partite sets $\{x_1,x_2\}$ and $\{y_1,y_2,y_3\}$. Let $D(8)$ be the bipartite digraph obtained from the digraph $K_{2,3}^*$ by adding  three new vertices $x_0,y_0, x_3$ and the following new arcs $x_0y_0$, $y_0x_0$, $x_0y_1$, $y_1x_0$, $x_3y_3$ and $y_3x_3$.

It is not difficult to check that the digraph $D(8)$ is strongly connected  and satisfies condition $B_0$, but the underlying undirected graph of $D(8)$ is not 2-connected  and $D(8)$ has no cycle of length 6.\\

\section {Proof of the main result}

\textbf{Proof of Theorem 1.6}. Let a digraph $D$ satisfies the conditions of the theorem. Suppose that  $D$ contains no cycle of length $2a-2$. Since the underlying undirected graph of $D$  is not 2-connected, it follows that  $V(D)=A\cup B\cup \{u\}$, where $A$ and $B$ are non-empty disjoint subsets of vertices of $D$, the vertex $u$ is not in $A\cup B$ and there is no arc between $A$ and $B$. Since $D$ is strong, there are vertices $x\in A$ and $x_0\in B$ 
such that $\{x,x_0\}\rightarrow u$, i.e., $\{x,x_0\}$ is a dominating pair. Note that $x$ and $x_0$ belong to the same partite set, say $X$. Then $u\in Y$. By  condition $B_0$,  $max\{d(x), d(x_0)\}\geq 2a-2$. 
Without loss of generality, we assume that $d(x)\geq 2a-2$. From this and the fact that there are no arc between $A$ and $B$ it follows that $a-2\leq |Y\cap A|\leq a-1$.

Put $Y_1:=Y\cap A$. We will consider the cases $|Y_1|=a-2$ and $|Y_1|=a-1$ separately.\\

\textbf{Case 1.}  $|Y_1|=a-2$.

Then $|Y\cap B|=1$. Let $Y_1:=\{y_1,y_2,\ldots , y_{a-2}\}$ and $Y\cap B:=\{y_0\}$. It is not difficult to check that the vertex $x$ and every vertex of $Y_1\cup \{u\}$ form a 2-cycle, i.e., $x\leftrightarrow Y_1\cup \{u\}$. Therefore every pair of distinct vertices of $Y_1\cup \{u\}$ is a dominating pair. This means that $Y_1\cup \{u\}$ has at least $a-2$ vertices (maybe except, say $y_{a-2}$, or $u$) each of which has degree at least $2a-2$. Then $d(y_1)\geq 2a-2$, since $a\geq 5$. From this it follows that $|X\cap A|=a-1$ and $X\cap B=\{x_0\}$ since there is no arc between $y_1$ and $B$. 

Put $X_1:=\{x_1,x_2,\ldots , x_{a-1}\}$, where $x_1=x$.
Therefore, $B=\{x_0,y_0\}$. Since $D$ is strong and  $y_0$ is not adjacent to any vertex of $X_1$, it follows that $y_0\leftrightarrow x_0$, $u\rightarrow x_0$, $d(x_0)=4$ and $d(y_0)=2$.
By condition $B_0$, $d(u)\geq 2a-2$ since $\{u,y_0\}\rightarrow x_0$.

First consider the case when $d(y_i)\geq 2a-2$ for all $i\in [1,a-2]$. Then $Y_1\leftrightarrow X_1$, since there is no arc between $Y_1$ and $\{x_0\}$, i.e., the induced subdigraph $ D\langle Y_1 \cap X_1 \rangle$ is a complete bipartite  digraph with partite sets $X_1$ and $Y_1$. Since $d(u)\geq 2a-2$, it follows that there are  at least $a-3$ vertices in $X_1$ each of which together with $u$ form a 2 cycle. Now we can choose a vertex in $X_1$ other than $x$, say $x_2$, such that $u\leftrightarrow x_2$. Therefore, $x_1ux_2y_2x_3\ldots x_{a-2}y_{a-2}x_{a-1}y_1x_1$ is a cycle of length $2a-2$, which contradicts the supposition that $D$ contains no cycle of length $2a-2$.

Now consider the case when there is   a vertex in $Y_1$, say $y_{a-2}$, which has degree at most $2a-3$. Then from condition $B_0$ it follows that $d(y_i)\geq 2a-2$ 
for all $i\in [1,a-3]$ since $x\leftrightarrow Y_1\cup \{u\}$ and $d(x)\geq 2a-2$. This implies that the subdigraph
 $ D\langle X_1 \cup \{y_1,y_2,\ldots , y_{a-3}\} \rangle$ is a complete bipartite digraph with partite sets $X_1$ and 
 $\{y_1,y_2,\ldots , y_{a-3}\}$. In particular, $y_1\leftrightarrow X_1$. Then every pair of distinct vertices of $X_1$ is a dominating pair. Condition $B_0$ implies  that  $X_1$ contains at least $a-2$ vertices, say $x_1,x_2, \ldots ,x_{a-2}$ , each of which has degree at least $2a-2$. 
Then 
$$
\{x_1,x_2,\ldots , x_{a-2}\}\leftrightarrow Y_1\cup \{u\},$$ in particular $y_{a-2}\leftrightarrow \{x_1,x_2,\ldots , x_{a-2}\}$ and 
$u\leftrightarrow \{x_1,x_2,\ldots , x_{a-2}\}$.
Therefore, $y_1x_{a-1}y_2x_2y_3x_3\ldots y_{a-2}$ $x_{a-2}ux_1y_1$ is a cycle of length $2a-2$, which is a contradiction.\\

\textbf{Case 2.}  $|Y_1|=a-1$.

Let now $Y_1:=\{y_1,y_2,\ldots , y_{a-1}\}$. Then $Y\cap B=\emptyset$, i.e., $B\subseteq X$. Since $D$ is strong, from condition $B_0$ it follows  that $B=\{x_0\}$, $u\leftrightarrow x_0$ and $|X\cap A|=a-1$. Let again $X_1:=X\cap A=
\{x_1,x_2,\ldots , x_{a-1}\}$, where $x_1=x$ (recall that $x_1\rightarrow u$).

If $d(y_i)\geq 2a-2$ for all $i\in [1,a-1]$,  then the subdigraph 
$D\langle X_1\cup Y_1\rangle$ is a complete bipartite digraph with partite sets $X_1$ and $Y_1$. Therefore, $D$ contains a cycle of length $2a-2$, a contradiction.
We may therefore assume that  $Y_1$ contains a vertex of degree at most $2a-3$. Observe that $Y_1$ may contains at most three vertices each of which has degree less than $2a-2$ since $d(x_1)\geq 2a-2$ (for otherwise $Y_1$ contains two vertices, say $v$ and $z$, such that $\{v,z\}\rightarrow x_1$ and $max\{d(v),d(z)\}\leq 2a-3$, which contradict condition $B_0$). We  consider the following three possible subcases
depending on the number of vertices in $Y_1$ each of which has degree at most $2a-3$.\\

\textbf{Subcase 2.1}. $Y_1$ contains exactly one vertex of degree less than $2a-2$.

Assume, without loss of generality, that $d(y_{a-1})\leq 2a-3$ and $d(y_i)\geq 2a-2$ for all $i\in [1,a-2]$. Then it is easy to see that the subdigraph $D\langle X_1\cup Y_1\setminus \{y_{a-1}\}\rangle$ is a complete bipartite digraph with partite sets $X_1$ and 
$Y_1\setminus \{y_{a-1}\}$ since $d(x_0, Y_1)=0$. From strong connectedness of $D$ it follows that $d^+(u,X_1)\geq 1$. 
If $u\rightarrow x_i$ for some $i\in [2,a-1]$, then   by the symmetry between of  vertices $x_2,x_3,\ldots , x_{a-1}$, we can assume that $u\rightarrow x_2$. Then it is easy to see that
$ux_2y_2x_3\ldots y_{a-2}x_{a-1}y_1x_1u$ is a cycle of length $2a-2$, which is a contradiction. We may therefore assume that
 $$
d^+(u,\{x_2,x_3,\ldots x_{a-1}\})=0.  \eqno (1)
$$
 Then $u\rightarrow x_1$,  $d^+(y_{a-1})\geq 1$ and $d^-(y_{a-1})\geq 1$, 
since $D$ is strong. If there exist two distinct vertices in $X_1$, say $x_1$ and $x_2$, such that
 $x_1\rightarrow y_{a-1}$ and  $y_{a-1}\rightarrow x_2$, then the cycle $x_1y_{a-1}x_2y_2x_3\ldots x_{a-2}y_{a-2}x_{a-1}y_1x_1$ is a cycle of length $2a-2$, a contradiction. We may therefore assume that there are no two distinct vertices $x_i$ and $x_j$ of $X_1$ such that $x_i\rightarrow y_{a-1}$ and  $y_{a-1}\rightarrow x_j$. Then $d^+(y_{a-1})=d^-(y_{a-1})= 1$ and $y_{a-1}\leftrightarrow x_i$ for some $i\in [1,a-1]$. If $i=1$, i.e., $x_1\leftrightarrow y_{a-1}$, then  $d(y_{a-1})=2$. Now using (1) and the fact that $d(u,\{x_0,x_1\})=4$, we obtain 
$$
d(u)=d(u,\{x_0,x_1\})+d^-(u,\{x_2,x_3,\ldots x_{a-1}\})\leq a+2\leq 2a-3,
$$
which contradicts condition $B_0$ since $\{u,y_{a-1}\}\rightarrow x_1$ and $a\geq 5$. Therefore $i\in [2,a-1]$.
  
Assume, without loss of generality, that $y_{a-1}\leftrightarrow x_{a-1}$. Then $a(x_i,y_{a-1})=0$ for all  $i\in [1,a-2]$, in particular, $a(x_2,y_{a-1})=a(x_3,y_{a-1})=0$.  
This together with 
(1) imply that $max\{d(x_2),d(x_3)\}\leq 2a-3$, which contradicts condition $B_0$ since
$\{x_2,x_3\}\rightarrow y_1$. The discussion of Subcase 2.1 is completed. \\  

 \textbf{Subcase 2.2}. $Y_1$ contains exactly two  vertices each of which has degree less than $2a-2$.

Assume, without loss of generality, that $d(y_{a-2})\leq 2a-3$, $d(y_{a-1})\leq 2a-3$ and $d(y_i)\geq 2a-2$ for all $i\in [1,a-3]$. Then it is easy to see that the subdigraph $D\langle X_1\cup Y_1\setminus \{y_{a-2},y_{a-1}\}\rangle$ is a complete bipartite digraph with partite sets $X_1$ and 
$Y_1\setminus \{y_{a-2},y_{a-1}\}$ since $d(x_0, Y_1)=0$.\\ 

 We prove the following  Claims 1 and 2 below.\\

\textbf{Claim 1}. If $x_j\rightarrow y_{a-2}$ for some $j\in [2,a-1]$, then $d^+(y_{a-2},\{x_1,x_2,\ldots , x_{a-1}\}\setminus \{x_j\})=0$.

\textbf{Proof of Claim 1}. Assume, without loss of generality, that $x_{a-1}\rightarrow y_{a-2}$, i.e., $j=a-1$. Suppose that the claim is not true, i.e., $y_{a-2}\rightarrow x_i$ for some $i\in [1,a-2]$. We will consider the cases $i=1$ and
$i\in [2,a-2]$ separately.

\textbf{Case}. $i=1$, i.e., $y_{a-2}\rightarrow x_1$.

First we show that 
$$
d^+(u,\{x_2,x_3,\ldots ,x_{a-1}\})=0.   \eqno (2)
$$

\textbf{Proof of (2)}. Suppose that (2) is not true, i.e., there is an $k\in [2,a-1]$ such that $u\rightarrow x_k$.  If $k\in [2,a-2]$, 
we may assume, without loss of generality, that $u\rightarrow x_2$. Then the cycle $x_{a-1}y_{a-2}x_1ux_2y_1x_3y_2\ldots $ $x_{a-2}y_{a-3}x_{a-1}$ is a cycle of length $2a-2$, contradiction. Thus, we may assume  that $k=a-1$. Then
$$
 u\rightarrow x_{a-1} \quad \hbox{and} \quad d^+(u,\{x_2,x_3,\ldots ,x_{a-2}\})=0 .  \eqno (3)
$$
If $y_{a-2}\rightarrow x_l$, for some $l\in [2,a-2]$ (say $y_{a-2}\rightarrow x_2$), then the cycle $x_{a-1}y_{a-2}x_2y_1x_3y_2\ldots x_{a-2}y_{a-3}x_1\\ux_{a-1}$ is a cycle of length $2a-2$, a contradiction. We may therefore assume  that 
$$
d^+(y_{a-2},\{x_2,x_3,\ldots ,x_{a-2}\})=0.   \eqno (4)
$$
If $x_l\rightarrow u$ for some $l\in [2,a-2]$ (say $x_2\rightarrow u$), then  the cycle 
$x_{a-1}y_{a-2}x_1y_2x_3y_3\ldots y_{a-3}x_{a-2}y_{1}x_2ux_{a-1}$ is a cycle of length $2a-2$, a contradiction.  We may therefore assume that 
 $$
d^-(u,\{x_2,x_3,\ldots ,x_{a-2}\})=0.   
$$
Combining this together with  (3) and (4), we obtain
$$
d(u,\{x_2,x_3,\ldots ,x_{a-2}\})=d^+(y_{a-2},\{x_2,x_3,\ldots ,x_{a-2}\})=0.
$$
This and  $a\geq 5$ imply that $d(x_2)\leq 2a-3$ and $d(x_3)\leq 2a-3$, which contradict condition $B_0$ since $\{x_2,x_3\}\rightarrow y_1$. This 
contradiction proves (2).\\

Since $D$ is strong, from (2) it follows that $u\rightarrow x_1$. Therefore, $\{u,y_{a-2}\}\rightarrow x_1$, 
i.e., $\{u,y_{a-2}\}$ is a dominating pair. Now using condition $B_0$, we obtain that  $d(u)\geq 2a-2$ since $d(y_{a-2})\leq 2a-3$ (by our assumption). Therefore, by (2), 
$$
2a-2\leq d(u)=d(u,\{x_0,x_1\})+d^-(u,\{x_2,x_3,\ldots ,x_{a-1}\})\leq 4+a-2=a+2.
$$
Hence, $a\leq 4$, which contradicts that $a\geq 5$. The discussion of the case $i=1$ is completed.\\

 \textbf{Case}.  $i\in [2,a-2]$, i.e., $y_{a-2}\rightarrow x_i$ and $y_{a-2}x_1\notin A(D)$.

Assume, without loss of generality, that $y_{a-2}\rightarrow x_2$, i.e., $i=2$. Now we prove that 
$$
d^+(u,\{x_3,x_4,\ldots ,x_{a-1}\})=0.   \eqno (5)
$$

\textbf{Proof of (5)}. Suppose that (5) is not true, i.e., there is an $l\in [3,a-1]$ such that  $u\rightarrow x_l$.  If $l=a-1$, i.e.,
$u\rightarrow x_{a-1}$, then the cycle $x_{a-1}y_{a-2}x_2y_2x_3\ldots y_{a-3}x_{a-2}y_{1}x_1ux_{a-1}$ is a cycle of length $2a-2$. We may  therefore assume that $l\in [3,a-2]$. Without loss of generality,  assume that $u\rightarrow x_3$. Then the cycle $x_1ux_3y_2x_4\ldots y_{a-4}x_{a-2}y_{a-3}$ $x_{a-1}y_{a-2}x_2y_1x_1$ is a cycle of length $2a-2$. In both cases we have a cycle of length $2a-2$, which is a contradiction. Therefore (5) is true.\\

From (5) and strongly connectedness of $D$ it follows that $u\rightarrow x_1$ or $u\rightarrow x_2$. 
 
First consider the case  $u\rightarrow x_1$.
It is not difficult to show that 
$$
d^-(u,\{x_2,x_3,\ldots ,x_{a-2}\})=0.   \eqno (6)
$$
Indeed, if $x_2\rightarrow u$, then the cycle $y_{a-2}x_2ux_1y_1x_3y_2x_4\ldots x_{a-2} y_{a-3}x_{a-1}y_{a-2}$ has length $2a-2$;
 if 
$x_j\rightarrow u$ and $j\in [3,a-2]$, then (we may assume that $j=3$, i.e., $x_3\rightarrow u$) the cycle 
$x_{a-1}y_{a-2}x_2y_1x_3ux_1y_2x_4\ldots $ $y_{a-4}x_{a-2}y_{a-3}x_{a-1}$ has length $2a-2$.  In both cases we have a contradiction. Therefore, the equality (6) is true.\\

It is not difficult to show that $ux_2\notin A(D)$. Assume that this is not the case, i.e., $ux_2\in A(D)$. Then   from 
 $y_{a-2}\rightarrow x_2$, $d(y_{a-2})\leq 2a-3$ and condition $B_0$ it follows that 
$d(u)\geq 2a-2$. On the other hand, using (5) and (6) we obtain
$$
2a-2\leq d(u)=d(u,\{x_0,x_1\})+d^+(u,\{x_2,x_3,\ldots ,x_{a-1}\})+d^-(u,\{x_2,x_3,\ldots ,x_{a-1}\})\leq 6.
$$
Therefore, $a\leq 4$, which contradicts that $a\geq 5$. Thus,  $u x_2\notin A(D)$. This together with  (5) and (6) imply that
$$
d^+(u,\{x_2,x_3,\ldots ,x_{a-1}\})=d^-(u,\{x_2,x_3,\ldots ,x_{a-2}\})=0.
$$
In particular, $a(x_j,u)=0$ for all $j\in [2,a-2]$. Since $a\geq 5$ and  $\{x_2,x_3\}\rightarrow y_1$, it follows 
that $d(x_2)=2a-2$ or $d(x_3)=2a-2$. If $d(x_2)=2a-2$, then  $\{y_{a-2},y_{a-1}\}\rightarrow x_2$, and if   $d(x_3)=2a-2$, then $\{y_{a-2},y_{a-1}\}\rightarrow x_3$. In each case we have a   contradiction to condition $B_0$ because of 
$d(y_{a-2}\leq 2a-3$ and $d(y_{a-1}\leq 2a-3$.\\

Consider now the case $u x_1\notin A(D)$ and $u\rightarrow x_2$. 
Then, by condition $B_0$,  $d(u)\geq 2a-2$ since $\{u,y_{a-2}\}\rightarrow x_2$ and $d(y_{a-2})\leq 2a-3$. Now using (5),
we obtain
$$
2a-2\leq d(u)=d(u,\{x_0,x_1\})+d^+(u,\{x_2,x_3,\ldots ,x_{a-1}\})+d^-(u,\{x_2,x_3,\ldots ,x_{a-1}\})\leq a+2,
$$
which is a contradiction, because of $a\geq 5$. Claim 1 is proved. \fbox \\\\

\textbf{Claim 2}. If $x_j\rightarrow y_{a-2}$ for some $j\in [2,a-1]$, then $d^-(y_{a-2},\{x_2,x_3,\ldots , x_{a-1}\}\setminus \{x_j\})=0$.

\textbf{Proof of Claim 2}. Assume, without loss of generality, that $x_{a-1}\rightarrow y_{a-2}$, i.e., $j=a-1$. Suppose that  the claim is not true, i.e., $x_l\rightarrow y_{a-2}$ for some $l\in [2,a-2]$. Assume, without loss of generality, that $x_{2}\rightarrow y_{a-2}$, i.e., $l=2$. From Claim 1 and strongly connectedness
of $D$ it follows that  $y_{a-2}\rightarrow x_{a-1}$.  This together with condition $B_0$ and 
$max\{d(y_{a-2}),d(y_{a-1})\}\leq 2a-3$ imply that $y_{a-1}x_{a-1}\notin A(D)$. 
 If $u\rightarrow x_2$, then the cycle
$y_{a-2}x_{a-1}y_2x_3y_3\ldots x_{a-3}y_{a-3}x_{a-2}y_1x_1ux_2y_{a-2}$ has length $2a-2$, which is a contradiction. If 
$u\rightarrow x_k$, where $k\in [3,a-2]$. We may assume that $k=3$, i.e., $u\rightarrow x_3$. Then 
$y_{a-2}x_{a-1}y_1x_1ux_3y_2x_4y_3\ldots x_{a-2}y_{a-3}$ $x_2y_{a-2}$ is a cycle of length $2a-2$, which is a contradiction.
We may therefore assume that 
$$
d^+(u,\{x_2,x_3,\ldots ,x_{a-2}\})=0.   \eqno (7)
$$
From (7) and strongly connectedness of $D$ it follows that $u\rightarrow x_1$ or $u\rightarrow x_{a-1}$.

First consider the case $u\rightarrow x_1$. It is not difficult to see that if for some $j\in [3,a-2]$, say $j=3$, 
$x_j\rightarrow u$, then the cycle 
$y_{a-2}x_{a-1}y_1x_3ux_1y_3x_4\ldots y_{a-3}x_{a-2}y_2x_2y_{a-2}$ has length $2a-2$, and if $x_{a-1}\rightarrow u$, then the cycle 
$y_{a-2}x_{a-1}ux_1y_1x_3y_3\ldots x_{a-3}y_{a-3}x_{a-2}y_2x_2y_{a-2}$ has length $2a-2$, which is a contradiction.
We may therefore assume that
$$
d^-(u,\{x_3,x_4,\ldots ,x_{a-1}\})=0.   \eqno (8)
$$
Now using (7) and (8), we obtain that $a(u,x_j)=0$ for all $j\in [3,a-2]$ and 
$$
d(u)=d(u,\{x_0,x_1\})+d^+(u,\{x_2,x_3,\ldots ,x_{a-1}\})+d^-(u,\{x_2,x_3,\ldots ,x_{a-1}\})\leq 6\leq 2a-3.
$$
From (7), (8) and Claim 1 it follows that $d(x_j)\leq 2a-3$ for all $j\in [3,a-2]$. Hence, $a-2=3$, i.e., $a=5$ and $d(x_3)\leq 2a-3$. By condition $B_0$ and $\{x_2,x_3,\ldots , x_{a-1}\}\rightarrow y_1$, we have that $d(x_2)\geq 2a-2$ and 
$d(x_4)\geq 2a-2$.  From $y_{a-1} x_{a-1}\notin A(D)$ and $x_{a-1} u\notin A(D)$ ($a-1=4$) it follows that $u\rightarrow x_{a-1}$, which is a contradiction since $\{u,y_{a-2}\}\rightarrow x_{a-1}$ and $max\{d(u), d(y_{a-2}\}\leq 2a-3$.

Now consider the case when $u\rightarrow x_{a-1}$ and $u x_1\notin A(D)$. Since $\{u, y_{a-2}\}\rightarrow x_{a-1}$ and 
$d(y_{a-2})\leq 2a-3$, it follows that $d(u)\geq 2a-2$. On the other hand, using (7) and $ux_1\notin A(D)$, we obtain
$$
2a-2\leq d(u)=d(u,\{x_0,x_1\})+d^+(u,\{x_2,x_3,\ldots ,x_{a-1}\})+d^-(u,\{x_2,x_3,\ldots ,x_{a-1}\})\leq a+2,
$$
which contradicts that $a\geq 5$. Claim 2 is proved. \fbox \\\\

Now we are ready to complete the discussion of Subcase 2.2.

Assume that $d^-(y_{j},\{x_2,x_3,\ldots , x_{a-1}\})\not= 0$ for $j=a-2$ or $j=a-1$ (say $j=a-2$). Assume, without loss of generality, that $x_{a-1}\rightarrow y_{a-2}$.
From Claims 1 and 2 it follows that 
$$
d^+(y_{a-2},\{x_1,x_2,x_3,\ldots ,x_{a-2}\})=d^-(y_{a-2},\{x_2,x_3,\ldots ,x_{a-2}\})=0. \eqno (9)
$$
Therefore, $d(x_i)\leq 2a-2$ for all $i\in [2,a-2]$ since $a(x_i,y_{a-2})=0$. From strongly connectedness of $D$ and (9) it follows that $y_{a-2}\rightarrow x_{a-1}$. This along with $max\{d(y_{a-2}), d(y_{a-1})\}\leq 2a-3$ and condition $B_0$ implies that
 $y_{a-1}x_{a-1}\notin A(D)$. Therefore, 
$$
d^+(y_{a-1},\{x_1,x_2,x_3,\ldots ,x_{a-2}\})\not=0
$$
 since $D$ is strong.  Applying Claim 1 to the vertex $y_{a-1}$
we obtain that $x_{a-1} y_{a-1}\notin A(D)$. Then $a(x_{a-1},y_{a-1})=0$ and $d(x_{a-1})\leq 2a-2$. Since
$\{x_2,x_3,\ldots ,x_{a-1}\} \rightarrow y_1$, from condition $B_0$ it follows that  
$\{x_2,x_3,\ldots ,x_{a-1}\}$ contains at least $a-3$ vertices each of which has degree at least $2a-2$. In particular, $d(x_2)\geq 2a-2$ or $d(x_3)\geq 2a-2$.
Without loss of generality, we assume that  $d(x_2)\geq 2a-2$. Then $x_2\rightarrow \{u,y_{a-1}\}$ since 
$a(x_2,y_{a-2})=0$. Now applying Claims 1 and 2 respect to the vertex   $y_{a-1}$, similarly to (9), we obtain 
$$
d^+(y_{a-1},\{x_1,x_3,x_4,\ldots , x_{a-1}\})=d^-(y_{a-1},\{x_3,x_4,\ldots ,x_{a-1}\})=0. \eqno (10)
$$
In particular,  from (9) and (10) we have $d^-(x_1,\{y_{a-2},y_{a-1}\})=0$. 
Therefore $x_1\rightarrow y_{a-2}$ and $u\rightarrow x_{1}$ because of $d(x_1)\geq 2a-2$. Hence, the cycle
$x_{2}ux_1y_{a-2}x_{a-1}y_1x_3y_2x_4\ldots y_{a-4}x_{a-2}y_{a-3}x_2$ has length equal to $2a-2$, which is a contradiction.

Now consider the case 
$$
A(\{x_2,x_3,\ldots , x_{a-1}\}\rightarrow \{y_{a-2},y_{a-1}\})= 0.  
$$ 
 Then, since $D$ is strong,  
 it follows that $x_1\rightarrow \{y_{a-2},y_{a-1}\}$. From the last equality we have that $d(x_j)\leq 2a-2$ for all
$j\in [2,a-1]$. This together with $\{x_2,x_3,\ldots ,x_{a-1}\} \rightarrow y_1$ imply that there are  at least $a-3$ vertices in
$\{x_2,x_3,\ldots ,x_{a-1}\}$ each of which has degree equal to $2a-2$. Assume, without loss of generality, that $d(x_2)=2a-2$. 
Then $\{y_{a-2},y_{a-1}\} \rightarrow x_2$, which is a contradiction since $d(y_{a-2})\leq 2a-3$ and $d(y_{a-1})\leq 2a-3$.
 In each case we obtain a contradiction, and hence the discussion of Subcase 2.2 is completed.\\

\textbf{Subcase 2.3}. $Y_1$ contains exactly three  vertices each of which has degree less than $2a-2$.

Assume, without loss of generality, that $d(y_{j})\leq 2a-3$ for all $j\in [a-3,a-1]$  and $d(y_i)\geq 2a-2$ for all $i\in [1,a-4]$.
 Then it is easy to see that the subdigraph $D\langle X_1\cup \{y_1,y_2,\ldots , y_{a-4}\}\rangle $ is a complete  bipartite digraph and $d^-(x_i,\{y_{a-3},y_{a-2},y_{a-1}\})\leq 1$ for all $i\in [1,a-1]$. This together with condition $B_0$ imply that  $\{x_2,x_3,\ldots ,x_{a-1}\}$ contains at least $a-3$ vertices, say $x_2,x_3,\ldots , x_{a-2}$, each of which has degree equal to $2a-2$. Then $x_1\leftrightarrow u$ (since $d(x_1)\geq 2a-2$ by our assumption), $x_i\rightarrow \{y_{a-3},y_{a-2},y_{a-1}\}$ if $i\in [1,a-2]$, and 
$x_j\leftrightarrow u$ if $j\in [2,a-2]$. Now it is not difficult to see that for every $i\in [1,a-2]$ there is an $j\in [a-3,a-1]$ such that $x_i\leftrightarrow y_j$. Because of the symmetry between of vertices $x_1,x_2,\ldots , x_{a-2}$ and the symmetry between of vertices $y_{a-3}$, $y_{a-2}$, $y_{a-1}$ we can assume,    
$x_1\leftrightarrow y_{a-3}$.

 First  consider the case when 
$$
A(\{y_{a-2},y_{a-1}\}\rightarrow \{x_4,x_5,\ldots ,x_{a-1}\})\not=\emptyset.
$$ 
By the symmetry between of vertices $x_4,x_5,\ldots , x_{a-1}$ and the symmetry between of vertices  $y_{a-2}$, $y_{a-1}$ we can assume that    
 $y_{a-2}\rightarrow x_{a-1}$. Therefore, if $a\geq 6$, then the cycle $x_2ux_3y_{a-3}$ $x_1y_{a-2}x_{a-1}y_1x_4y_2\ldots $ $ x_{a-2}y_{a-4}x_2$ has length $2a-2$, and if $a=5$, then the cycle 
$x_2ux_3y_{2}x_1y_3x_4y_1x_2$ has length $2a-2$, which is a contradiction.

Now consider the case when
 $$
A(\{y_{a-2},y_{a-1}\}\rightarrow \{x_4,x_5,\ldots ,x_{a-1}\})=\emptyset.  
$$
Then, since $D$ is strong, from $x_1\leftrightarrow y_{a-3}$, 
$max\{d(y_{a-3}),d(y_{a-2}),d(y_{a-1})\}\leq 2a-3$ and condition $B_0$  it follows that 
$$
d^-(x_1,\{y_{a-2},y_{a-3}\})=0 \quad \hbox{and} \quad 
min\{d^+(y_{a-2},\{x_2,x_3\}),d^+(y_{a-1},\{x_2,x_3\}\})\geq 1. \eqno (11)
$$
Without loss of generality, we assume that $y_{a-2}\rightarrow x_2$. If $a\geq 6$, then the cycle 
$y_{a-2}x_2y_{a-3}x_1ux_3y_1x_{a-1}$ $y_2x_4y_3\ldots x_{a-3}y_{a-4}x_{a-2}y_{a-2}$ has length $2a-2$, which is a contradiction. We may therefore assume  that $a=5$. By the above observation we have  that $D\langle \{x_1,x_2,x_3,u,y_1\}\rangle$ is a complete bipartite digraph with partite sets $\{x_1,x_2,x_3\}$ and $\{u,y_1\}$, $\{x_1,x_2,x_3\}\rightarrow \{y_2,y_3,y_4\}$,
 $x_4\leftrightarrow y_1$, $x_0\leftrightarrow u$ and $x_i\leftrightarrow y_{i+1}$ for all $i\in [1,3]$. From $d(y_j)\leq 2a-3$, where $j\in [2,4]$, we obtain that
$$
d^-(x_1,\{y_3,y_4\})=d^-(x_2,\{y_2,y_4\})=d^-(x_3,\{y_2,y_3\})=0.
$$
Using these,  it is not difficult to show that $d(x_4,\{u,y_2,y_3,y_4\})=0$ (for otherwise, $D$ would contain a cycle of length eight, which contradicts our initial supposition). Now
it is not difficult to check that the obtained digraph is strongly connected and isomorphic to $D(10)$, which satisfies condition $B_0$, but has no cycle of length 8.
 \fbox \\\\

From Theorems 1.5 and 1.6 it follows the following corollary. 

\textbf{Corollary}. {\it Let $D$ be a strongly connected balanced bipartite digraph of order $2a\geq 8$. Assume that the underlying undirected graph of $D$  is not 2-connected and $d(x)+d(y)\geq 4a-3$ for every dominating pair of vertices $x$ and $y$.  Then  $D$ contains a cycle of length $k$ for every $k\in [1,a-1]$  unless $D$ is isomorphic to the digraph $D(10)$.}\\

\end{document}